\documentclass[letterpaper, 10 pt, conference]{ieeeconf}  % Comment this line out
\IEEEoverridecommandlockouts                              % This command is only
\usepackage{graphics} % for pdf, bitmapped graphics files
\usepackage{epsfig} % for postscript graphics files

\usepackage{amssymb,amsfonts,amsmath}
\usepackage{latexsym}
\usepackage{mathrsfs}

\newcommand{\be}{\begin{equation}}
\newcommand{\ee}{\end{equation}}
\newcommand{\bea}{\begin{eqnarray}}
\newcommand{\eea}{\end{eqnarray}}
\newcommand{\beas}{\begin{eqnarray*}}
\newcommand{\eeas}{\end{eqnarray*}}
\newcommand{\pa}{\partial}
\newcommand{\parsh}[2]{\frac{\pa #1}{\pa #2}}
\newcommand{\la}{\lambda}

\newcommand{\ra}{\rightarrow}

\newtheorem{theorem}{Theorem}

\newtheorem{corollary}{Corollary}
\newtheorem{definition}{Definition}
\newcommand{\e}{\mathbf{e}}
\usepackage{psfrag}

\title{\LARGE \bf
Dynamics of Symplectic SubVolumes
}

\author{J.M. Maruskin${}^*$, D.J. Scheeres${}^\dag$, and A.M. Bloch${}^\ddag$
	\thanks{$*$ Ph.D. Candidate in Applied and Interdisciplinary Mathematics, University of Michigan.}
	\thanks{$\dag$ Associate Professor of Aerospace Engineering, University of Michigan}
	\thanks{$\ddag$ Alexander Ziwet Collegiate Professor of Mathematics and Department Chair, University of Michigan}}

\begin{document}

\maketitle
\thispagestyle{empty}
\pagestyle{empty}

%%%%%%%%%%%%%%%%%%%%%%%%%%%%%%%%%%%%%%%%%%%%%%%%%%%%%%%%%%%%%%%%%%%%%%%%%%%%%%%%
\begin{abstract}

In this paper we will explore fundamental constraints on the evolution of certain symplectic subvolumes
possessed by any Hamiltonian phase space.  This research has direct application to optimal control and control
of conservative mechanical systems.  We relate geometric invariants of symplectic topology to computations
that can easily be carried out with the state transition matrix of the flow map.  We will show how
certain symplectic subvolumes have a minimal obtainable volume; further if the subvolume dimension
equals the phase space dimension, this constraint reduces to Liouville's Theorem.  Finally we present
a preferred basis that, for a given canonical transformation, has certain minimality properties 
with regards to the local volume expansion of phase space.  
\end{abstract}

\section{Introduction}

\subsection{Overview}

The traditional approach for studying the dynamics and control of mechanical systems is to focus
on individual trajectories and states in order to determine where they will go and where they can be forced to
go.  In reality, however, system states are never precisely known and can only be determined to exist
within some set of finite volume in the dynamical system's phase space.  By treating such systems 
as a sum of individual trajectories, one loses the geometrical insight and deeper results offered by more wholistic
approaches.

In this paper, we will be concerned with understanding fundamental constraints on
the evolution of compact $2k$-dimensional
symplectic sets that evolve along a nominal trajectory of the system.  
Different symplectic constraints arise on such sets, including conservation
of the signed $2k$-volume projections on the coupled symplectic planes as well as 
the constraints implied from Gromov's Nonsqueezing Theorem (see Scheeres et al \cite{scheeres03} for a discussion
of these constraints in relation to orbit uncertainty evolution).
We will further present an additional constraint for a minimal obtainable volume 
that exists on certain classes of $2k$-dimensional symplectic sets and show how such a constraint
leads to the local collapse of phase space along solution curves in Hamiltonian phase space.
This collapse of phase space is fundamentally linked to the expansion of symplectic subvolumes.
Finally, for any fixed final time, no matter how large, we will produce a distinguished orthogonal symplectic
basis that resists collapse.  The basis may collapse as time evolves, but will return to being orthogonal at the final time.
The uncertainty of any $2k$-dimensional distributions initially parallel to the symplectic planes of this basis,
even though it may increase dramatically during the course of its evolution, will always return to its
initial uncertainty at the final time.

Since the resulting equations produced by applying Pontryagin's Maximum Principle to optimal control problems
are Hamiltonian, the results we discuss here should provide geometric insight to the evolution and control of uncertainty 
distributions in such systems.  
This theory provides fundamental limits on dynamical orbits, and hence
if one provides a control it provides limits on the accuracy of the control
in the face of uncertainty.  It also provides limits
on uncertainty propagation in optimal control systems.  
Moreover, the preferred minimal uncertainty basis we produce should have numerous benefits to the design of 
fixed finite time optimal control laws where precise state information is unknown.

\subsection{Outline}

In Section 2, we introduce Hamiltonian systems and the state transition matrix (STM).  We show how classical identities on the Lagrange and Poisson
brackets relate to constraints on the STM.  Specifically, for any symplectic column of the STM $\Phi$, the sum of the $2\times 2$ symplectic subdeterminants
must add up to unity.  

In Section 3, we look in depth at surfaces that can be explicitly parameterized by one of their symplectic planes.  We derive area
expansion factors from the parameterization plane to the surface, its image under the Hamiltonian phase flow, and the symplectic
projections of its image.  If the state of the system is somewhere on the initial surface (with equal a priori probability), we
interpret these various expansion factors as a probability map that leads one to understand where the particle, after applying the Hamiltonian
phase flow, is likely to be found.  

In Section 4, we present a nonsqueezing-like property for $2k$-volumes which is closely related to Wirtinger's Inequality.  We show how this leads to the
fact that subvolume expansions in the differential neighborhood of the Hamiltonian flow leads to a collapsing property in systems
which exhibit chaos.

In Section 5 we will discuss how these constraints on subvolume expansions, when considered with Liouville's Theorem, leads to 
the local collapse of the phase space around nominal trajectories.  Interestingly, we  will also show, 
that given a canonical transformation, there exists a preferred basis that resists collapse.  In particular, the volume of a $2k$ subvolume
chosen to be initially parallel to $k$ of the symplectic planes will return to its initial value
at this fixed final time.

\section{Hamiltonian Systems}

\subsection{The Classical Approach}

\subsubsection*{Hamilton's Equations}

In an $N$ degree of freedom Hamiltonian System, one has a $2N$-dimensional phase space spanned by $N$ generalized coordinates
$\{q_i\}_1^N$ and their $N$ conjugate momenta $\{ p_i \}_1^N$.  The conjugate pairings of coordinates and momenta, $(q_i, p_i)$
form what are known as symplectic pairs.  The dynamical equations of motion are deriveable from a Hamiltonian function $H(q,p)$
and Hamilton's equations:
\be
\label{sv01}
\dot q_i = \parsh{H}{p_i} \qquad \mbox{and} \qquad
\dot p_i = -\parsh{H}{q_i}
\ee
Solutions curves of the system \eqref{sv01} are called the Hamiltonian phase flow, and are denoted
$\phi_t(q, p)$.

\subsubsection*{Matrix Formalism}

Let $x = \langle p_1, q_1, p_2, q_2, \ldots p_N, q_N \rangle$ be the phase space position of the system, ordered by symplectic pairs.
Define the matrix
\[
J_2 = \left( \begin{array}{cc}
0 & -1 \\ 1 & 0 \end{array} \right)
\]
The symplectic matrix $J$ is defined as the $2N \times 2N$ block-diagonal matrix with $J_2$'s down the main diagonal.  Grouping the coordinates in 
symplectic pairs will be useful to us later as we will be looking at various ``symplectic columns'' of the STM.  Given the above definitions,
Hamilton's equations
can be cast into the following matrix form:
\be
\dot x = \frac{d}{dt} \phi_t(x_0) = J \cdot \parsh{H}{x}
\ee
where $x = \phi_t(x_0)$.
The right hand side is the so-called \textit{symplectic gradient} of the Hamiltonian function.

\subsubsection*{Lagrange Brackets}

Given a transformation 
\[
Q_i = Q_i(q, p) \qquad \mbox{and} \qquad P_i = P_i(q,p)
\]
we may introduce the \textit{Lagrange bracket} expression for the two variables $(u,v)$, 
(which can take on any of the values $q_1,  \ldots, q_n, p_1, \ldots, p_n$):
\be
\label{sv02}
[u, v] = \sum_{i=1}^N \left( \parsh{P_i}{u} \parsh{Q_i}{v} - \parsh{Q_i}{u} \parsh{P_i}{v} \right),
\ee
The exactness conditions required for a canonical transformation can then be cast into the following equivalent conditions:
\be
\label{sv03}
[q_j, q_k] = 0 \qquad [p_j, p_k] = 0 \qquad [p_j, q_k] = \delta_{jk}
\ee
The Hamiltonian phase flow $\phi_t(q, p)$ is a continuous one parameter family of canonical transformations.

\subsubsection*{Poisson Brackets}

Alternatively, one can define the Poisson bracket as:
\be
\label{sv04}
\{ u,v\} = \sum_{i=1}^N \left( \parsh{u}{p_i} \parsh{v}{q_i} - \parsh{u}{q_i} \parsh{v}{p_i} \right)
\ee
where $(u, v)$ can now be any of the variables $Q_1, \ldots, Q_n, P_1, \ldots, P_n$.  The sufficient conditions
for a canonical transformation can also be written as follows, in terms of the Poisson bracket:
\[
\{ Q_j, Q_k \} = 0 \qquad \{ P_j, P_k \} = 0 \qquad \{ P_j, Q_k \} = \delta_{jk}
\]

\subsection{The Geometric Approach}

\subsubsection*{Symplectic Manifolds}

A \textit{symplectic structure} on an even-dimensional manifold $M$ is a closed nondegenerate
differential two-form $\omega$ on $M$:
\[
d\omega = 0 \qquad \mbox{and} \qquad \forall \xi \not=0, \ \exists \eta : \omega(\xi, \eta) \not = 0
\]
The form $\omega$ is called the symplectic form and the pair $(M, \omega)$ is called a symplectic manifold.

On any symplectic manifold, there exists a vector space isomorphism between its cotangent and tangent bundles.
At $x \in M$, we have
\[
I_x: T^*_x  M \ra T_xM
\]
defined by the symplectic form and the following relation.  A vector $\xi \in T_xM$ is mapped to the 
one-form $I^{-1}_x(\xi)$ which acts on a vector $\eta \in T_xM$ as follows: $I^{-1}_x(\xi)(\eta) = \omega(\eta, \xi)$.

\subsubsection*{Hamiltonian Flows}

Let $H$ be a function $H: M \ra \mathbb{R}$ which we will call the Hamiltonian.  The associated Hamiltonian vector field
on $M$ is defined by $IdH$.  The flow generated by the vector field $I dH$ is the Hamilton phase flow $\phi_t$.  
If $M = \mathbb{R}^{2N}$ with the standard symplectic form $\omega_0 = \sum_{i}  dp_i \wedge dq_i$, we recover
Hamilton's equations \eqref{sv01}.

A transformation $\phi:M \ra M$ is considered \textit{symplectic} or \textit{canonical} if it preserves the symplectic form,
i.e., $\phi^* \omega = \omega$.  The Hamiltonian phase flow $\phi_t$ is a one parameter family of canonical transformations.

\subsubsection*{Integral Invariants}

A differential $k$-form $\alpha$ is an \textit{integral invariant} of the map $\phi$ if the integrals of $\alpha$ on any $k$-chain $\sigma$
is preserved as follows:
\[
\int_{\phi(\sigma)} \alpha = \int_\sigma \alpha
\]
The symplectic form $\omega$ is an integral invariant of the Hamiltonian flow.

To gain a physical intuition for what the integral invariant $\int \omega$ represents,
consider now a closed parametrized surface $\psi(\sigma)$ in $\mathbb{R}^{2n} = 
(p,q)$, with a parametrization given by 
$\psi: \sigma \subset \mathbb{R}^2 \longrightarrow \mathbb{R}^{2n}$, 
$\psi: (u,v) \longrightarrow
(q(u,v), p(u,v))$.  Then
\[
\int \!\!\! \int_{\psi(\sigma)} dp \wedge dq = \sum_{i=1}^n \int \!\!\! 
\int_\sigma  \parsh{(p_i, q_i)}{(u,v)} \ du \ dv
\]

$\int \!\!\! \int_{\sigma} \parsh{(p_i,q_i)}{(u,v)} du dv$ represents the oriented area of the projection
of the surface $\psi(\sigma)$ on the $i$-th symplectic plane.  By considering $\psi(\sigma)$ as an initial
surface and applying the Hamiltonian phase flow, mapping the surface to $\phi_t(\psi(\sigma))$; we recognize,
as a physical interpretation of the preservation of the symplectic form under canonical mappings, 
that the sum of the oriented areas of the projections onto the $N$ symplectic planes is preserved.

\subsection{The State Transition Matrix}

\subsubsection*{Definition}

If $\phi:M \ra M, \phi(p, q) = (P, Q)$ is a canonical transformation, its differential 
\[
d\phi: T_{(p,q)}M \ra T_{(P,Q)} M
\]
is, when represented in matrix form, known as the state transition matrix (STM) $\Phi$, a terminology adopted
from Linear Systems Theory.
Supposing $\phi_t$ sends $\phi_t: x_0 \ra \phi(x_0)$, the STM $\Phi$ maps initial deviations in the 
initial conditions $\Phi: \delta x_0 \ra \Phi \cdot \delta x_0$ to its final state, so that
$\phi_t(x_0 + \delta x_0) \approx x + \Phi \cdot \delta x_0$, to first order.  Here $x = (p, q) \in M$.

\subsubsection*{Dynamics}

If $\phi_t$ is the Hamiltonian phase flow, we have:
\[
\frac{d}{dt} \phi_t(x_0) = J \cdot \parsh{H}{x}(x)
\]
If we perturb the initial conditions to $x_0 + \delta x_0$, we find:
\[
\frac{d}{dt} \phi_t(x_0+\delta x_0) = J \cdot \parsh{H}{x}(x + \delta x)
\]
By expanding this in a Taylor Series, one sees:
\be
\frac{d}{dt} \Phi = J \cdot \parsh{^2 H}{x^2} \cdot \Phi
\ee
This defines a system of $4N^2$ differential equations that can be integrated numerically, simultaneously along with the nominal
solution curve $\phi_t(x_0)$.

\subsubsection*{Relation to Lagrange and Poisson Brackets}

We will relate the Lagrange and Poisson Brackets to determinants of various submatrices of the STM.  We will arrange the coordinates
in a symplectic order, so that $x = \langle p_1, q_1, \ldots, p_N, q_N \rangle$ and $X = \langle P_1, Q_1, \ldots, P_n, Q_N \rangle$.
In this fashion, the STM $\Phi$ is thought of as
\[
\Phi = \parsh{X}{x}
\]
Define the following subdeterminants:
\[
M_{ij} = \mathrm{det} \left( \begin{array}{cc}
\parsh{P_i}{p_j} & \parsh{P_i}{q_j} \\
\parsh{Q_i}{p_j} & \parsh{Q_i}{q_j}
\end{array} \right) = \parsh{P_i}{p_j} \parsh{Q_i}{q_j} - \parsh{P_i}{q_j} \parsh{Q_i}{p_j}
\]
Hence, $M_{ij}$ is the subdeterminant of the intersection 
of the $i$th \textit{symplectic} row of the STM with its $j$th \textit{symplectic} column.

It is easy to see the Lagrange and Poisson brackets are related to these subdeterminants as follows:
\bea
[p_j, q_j] &=& \sum_{i=1}^N M_{ij} \\
(P_i, Q_i) &=& \sum_{j=1}^N M_{ij}
\eea

\section{Symplectic Surfaces}

\subsection{Surface Classifications}

We will begin by making the following fairly natural definitions.

\begin{definition} A \textbf{globally symplectic surface} is a two-dimensional
submanifold of the phase space $\mathbb{R}^{2n}$ which admits a 1-to-1 mapping
to at least one of the symplectic planes via the projection operator, i.e.
it is a surface which can be parameterized in explicit form by one of its
symplectic coordinate pairs.
\end{definition}

This characteristic is not an invariant one.  It is possible, for example, for a lamina parallel to a symplectic plane to fold under some symplectic map,
so that its image under the map is not 1-to-1 with any symplectic plane.  

Another surface type we will consider is the following:

\begin{definition}
A \textbf{parametrically symplectic surface} (or \textbf{parasymplectic surface}, for short) 
is a 2-dimensional submanifold of $\mathbb{R}^{2n}$ that admits a parameterization that is a symplectic one, i.e.
one with a parameterization map that is canonical.
\end{definition}

The parasymplecticity of a surface is an invariant characteristic.  Let $\sigma$ be a lamina on the symplectic plane $(u,v)$
which is the parameterization of the surface $\phi(\sigma) \subset \mathbb{R}^{2n}$, where $\phi$ is the parameterization map.  
Let $\phi$ be a symplectomorphism, which exists if $\phi(\sigma)$ is a parasymplectic surface.  Let $\psi: \mathbb{R}^{2n} \ra \mathbb{R}^{2n}$
be an arbitrary symplectomorphism which takes $\psi: \phi(\sigma) \ra \psi(\phi(\sigma))$.  Then $\psi(\phi(\sigma))$ is parasymplectic,
with symplectic parameterization $\psi \circ \phi : \sigma \ra \psi(\phi(\sigma))$.

We will consider $\pi_i : \mathbb{R}^{2n} \ra \mathbb{R}^2$ the $i$-th \textit{symplectic} projection operator, so that
$\pi_i(\langle p_1, q_1, \ldots, p_N, q_N \rangle) = \langle p_i, q_i \rangle$.

\subsection{Area Expansion Factors}

\subsubsection*{Notation}

We define the $2n \times 2$ matrix $\Pi_\kappa$ as:
\[
\Pi_\kappa = \left[ \begin{array}{cccccccc}
\mathbb{O}_2 & \mathbb{O}_2 & \cdots & \mathbb{O}_2 & \mathbb{I}_2 & \mathbb{O}_2 & \cdots & \mathbb{O}_2
\end{array} \right]^T
\]
where the $\mathbb{I}_2$ appears in the $\kappa$-th symplectic row.
For any $2n \times 2n$ matrix $A$, the product $A \cdot \Pi_\kappa$ is the $\kappa$-th symplectic column
of $A$; $\Pi_\kappa^T \cdot A$ is the $\kappa$-th symplectic row of $A$; and $\Pi_\kappa^T \cdot A \cdot \Pi_\lambda$
is the $2 \times 2$ intersection of the $\kappa$-th symplectic row with the $\la$-th symplectic column.

In this notation, the STM subdeterminant $M_{ij}$, defined previously, can be expressed as:
\[
M_{ij} = \det(\Pi_i^T \cdot \Phi \cdot \Pi_j)
\]

\subsubsection*{Globally Symplectic Surfaces}

\begin{figure}[t]
\begin{center}
\psfrag{f}{$\phi_t$} \psfrag{dt}{$d\tau$}
\psfrag{pdt}{$\pi_j(d\tau)$} 
\psfrag{fdt}{$\phi_t(d\tau)$}
\psfrag{pfdt}{$\pi_i(\phi_t(d\tau))$}
\psfrag{pi}{$P_i$}
\psfrag{qi}{$Q_i$} 
\psfrag{pj}{$p_j$}
\psfrag{qj}{$q_j$}
\includegraphics[width=3in]{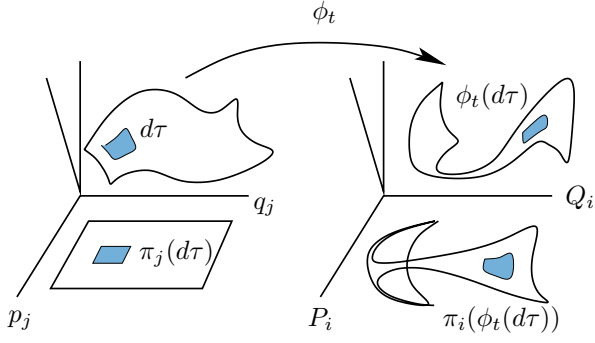}
\end{center}
\caption{Area Expansion Factors}
\label{svcdc1cdc1/fig01}
\end{figure}

We will consider now a surface $\tau$ which is \textit{globally symplectic} with respect to the $j$-th symplectic plane; i.e.,
the projection map $\pi_j: \tau \subset \mathbb{R}^{2n} \ra \pi_j(\tau) \subset \mathbb{R}^2$ is one-to-one.  We can parameterize
$\tau$ by its symplectic shadow on the $j$-th symplectic plane.  Now let the surface $\tau$ be mapped into the future by the Hamiltonian
flow $\phi_t: \langle p_1, q_1, \ldots, p_N, q_N \rangle \ra \langle P_1, Q_1, \ldots, P_N, Q_N \rangle$.  We will now consider the projection
of $\phi_t(\tau)$ onto the $i$-th symplectic plane.  For our analysis we will consider a differential area element $d\tau$ of $\tau$ 
(see Fig. \ref{svcdc1cdc1/fig01}).

The surface $\tau$ is described by the  parameterization $u=p_j$ and $v=q_j$ by $\langle p_1(u,v), q_1(u,v), \ldots,
u, v, \ldots, p_N(u,v), q_N(u,v) \rangle$.  
We define the matrix $L = \left[ \parsh{x_i}{u_j} \right] $ as follows:
\[
L = 
 \left[ \begin{array}{cc}
\parsh{p_1}{u} & \parsh{p_1}{v} \\
\parsh{q_1}{u} & \parsh{q_1}{v} \\
\vdots & \vdots \\
1 & 0 \\
0 & 1 \\
\vdots & \vdots \\
\parsh{p_N}{u} & \parsh{p_N}{v} \\
\parsh{q_N}{u} & \parsh{q_N}{v} \end{array} \right]
\]
where the $j$-th symplectic row is equated to $\mathbb{I}_2$.  $L$ is the Jacobian matrix of the parameterization
map that takes $\pi_j(d\tau) \ra d\tau$.  The metric for the surface $\tau$ in terms of variations in the $u-v$ plane
is given by:
\[
[g_{ij}] = \left[ \begin{array}{cc}
\parsh{\mathbf{x}}{u} \cdot \parsh{\mathbf{x}}{u} & \parsh{\mathbf{x}}{u} \cdot \parsh{\mathbf{x}}{v} \\
\parsh{\mathbf{x}}{v} \cdot \parsh{\mathbf{x}}{u} & \parsh{\mathbf{x}}{v} \cdot \parsh{\mathbf{x}}{v} \end{array} \right]
= L^T \cdot L
\]
where $\mathbf{x} = \langle p_1(u,v), \ldots, q_N(u,v) \rangle$, as usual.  Hence the metric determinant is:
\[
g=\det(g_{ij}) = \mathfrak{G}(L) = \det(L^T \cdot L)
\]
where $\mathfrak{G}(L)$ is the Gram determinant of the matrix $L$, defined by this equation.  

The surface area of $\tau$ is thus given by:
\[
S.A. = \int \!\!\! \int_\sigma \sqrt{g} du dv
\]

A simple application of this result gives us the total physical area of the area element $d\tau$ in terms of the area of its projection:
\[
\frac{\mathcal{A}(d\tau)}{\mathcal{A}(\pi_j(d\tau))} = \sqrt{\mathfrak{G}(L)}
\]
Similarly we have:
\[
\frac{\mathcal{A}(\phi_t(d\tau))}{\mathcal{A}(\pi_j(d\tau))} = \sqrt{\mathfrak{G}(\Phi \cdot L)}
\]
where $\Phi$ is the STM associated with $\phi_t$ and $x_0$.

The area of the projection $\pi_i(\phi_t(d\tau))$ is given by the Jacobian:
\[
\parsh{(P_i, Q_i)}{(u,v)} = \sum_{\kappa=1}^{2N} \sum_{\la=1}^{2N} \left(
\parsh{P_i}{x_\kappa} \parsh{Q_i}{x_\la} - \parsh{Q_i}{x_\kappa} \parsh{P_i}{x_\la} \right)
\parsh{x_\kappa}{u} \parsh{x_\la}{v}
\]
which can be represented more concisely as:
\[
\frac{\mathcal{A}(\pi_i(\phi_t(d\tau)))}{\mathcal{A}(\pi_j(d\tau))} = \det(\Pi_i^T \cdot \Phi \cdot L)
\]

\subsubsection*{Parasymplectic Surfaces}

We will now consider the case where $\tau$ is parallel to the $j$-th symplectic plane.  In this case, its parameterization map
is a symplectic one, and thus it is a parasymplectic surface.  All of the above results hold, but the matrix $L$ reduces to the 
simpler form $L = \Pi_j$, which gives us the following:
\[
\mathcal{A}(d\tau) = \mathcal{A}(\pi_j(d\tau))
\]
\[
\frac{\mathcal{A}(\phi_t(d\tau))}{\mathcal{A}(\pi_j(d\tau))} = \sqrt{\mathfrak{G}(\Phi \cdot \Pi_j)}
\]
and, most notably
\[
\frac{\mathcal{A}(\pi_i(\phi_t(d\tau)))}{\mathcal{A}(\pi_j(d\tau))} = \det(\Pi_i^T \cdot \Phi \cdot \Pi_j) = M_{ij}
\]
Preservation of the sum of the oriented symplectic area projections thus gives us the following constraint on the STM:
\be
\label{sv05}
\sum_{i=1}^N M_{ij} = [p_j, q_j] = 1
\ee

\subsubsection*{Application to Orbit Uncertainty Distributions}

Suppose we know that a system can be found anywhere on the surface $\tau$ with equal a priori probability.  The surface
is now mapped into the future by the Hamiltonian flow $\phi_t$ and we wish to determine where the particle is most likely to 
be on the $P_i-Q_i$ plane.  

We begin by discretizing the $(u=p_j)-(v=q_j)$ plane, each area element with area $\Delta u \Delta v$.  Summation will be assumed 
to be over each district.  The probability that the particle is in $d\tau$ is given by:
\[
\mathcal{P}(d\tau) = \frac{\mathcal{A}(d\tau)}{\mathcal{A}(\tau)} \approx 
\frac{\sqrt{\mathfrak{G}(L)} \Delta u \Delta v}{\sum \sqrt{\mathfrak{G}(L)} \Delta u \Delta v} = 
\frac{\sqrt{\mathfrak{G}(L)}}{\sum \sqrt{\mathfrak{G}(L)}}
\]
We reiterate the area of $\pi_i(\phi_t(d\tau))$ is given by:
\[
\mathcal{A}(\pi_i(\phi_t(d\tau))) \approx \det(\Pi_i^T \cdot \Phi \cdot L) \Delta u \Delta v
\]
so that the area probability density at $\langle P_i, Q_i \rangle$ is:
\[
\sigma \approx \frac{\sqrt{\mathfrak{G}(L)}}{|\det(\Pi_i^T \cdot \Phi \cdot L)| \Delta u \Delta v \sum \sqrt{\mathfrak{G}(L)}}
\]

This approach may be helpful in asteroid tracking, where angular and angular rate information is precisely known,
but there is initial uncertainty in the $r, \ \dot r$ distribution.  This problem is treated in Milani, et al. \cite{milani}.

\section{Integral Invariants}

In this section we will discuss the difference between two fundamental integral invariants defined for an arbitrary $2k$-dimensional
subvolume of our $2n$-dimensional phase space, for $k=1, \ldots, n$.  The first is the well-known integral invariant of Poincar\'e-Cartan.  
The second integral invariant is closely related, and is tantamount to a global version of Wirtinger's Inequality for lower dimensional 
subvolumes of phase space.

\subsection{Signed and Unsigned Integrals of Differential Forms}

We will begin our discussion on integral invariants with a brief discussion of 
the theory of integration of differential forms.  Let $\Sigma \subset \mathbb{R}^{2n}$ be a $2k$-dimensional
submanifold, parameterized by $\phi: (\sigma \subset \mathbb{R}^{2k}) \ra \Sigma$.  Let $\alpha$ be a $2k$-form on
$\mathbb{R}^{2n}$.  Then we define
\[
\int_{\Sigma} \alpha = \int_\sigma \phi^* \alpha
\]
For some function $f(x_1, \ldots, x_{2k})$, the pullback of $\alpha$ can be expressed in the following form:
\[
\phi^* \alpha = f(x_1, \ldots, x_{2k}) dx^1 \wedge \cdots \wedge dx^{2k}
\]
where we take $(x_1, \ldots, x_{2k})$ to be a basis of $\mathbf{R}^{2k}$.
Then the integral of $\phi^* \alpha$ over $\sigma$ reduces to the ordinary euclidean integral:
\[
\int_\Sigma \alpha = \int_{\sigma} \phi^* \alpha = \int_{\sigma} f(x_1, \ldots, x_{2k}) dx^1 \cdots dx^{2k}
\]

We would like to introduce a further definition as follows.  We define the unsigned integral of $\alpha$ over $\Sigma$
to be:
\[
\int_\Sigma |\alpha | = \int_{\sigma} |f(x_1, \ldots, x_{2k})| dx^1 \cdots dx^{2k}
\]
where $f$ has been defined above.

A more rigorous definition of the above integrals must involve a partition of unity, but simplicity has been choosen over rigor
so as to illustrate the spirit of the definitions.

\subsection{The Integral Invariants of Poincar\'e-Cartan}

Consider the standard symplectic form
\[
\omega = \sum_{i=1}^n p_i \wedge q_i
\]
and its $k$-th exterior product:
\[
\frac{1}{k!} \omega^k = \sum_{1 \le i_1 < \cdots < i_k \le n} dp_{i_1} \wedge dq_{i_1} \wedge \cdots \wedge dp_{i_k} \wedge dq_{i_k}
\]
Consider a set of $2k$ vectors $(X^1, \ldots, X^{2k})$ in $\mathbb{R}^{2n}$.  Then
\[
\frac{1}{k!} \omega^k(X^1, \ldots, X^{2k})
\]
represents the sum of the oriented $2k$-volume projections of the parallelpiped spanned by $X^1, \ldots, X^{2k}$ on
the symplectic ``$2k$-planes.''

$\omega^k$ is known as the integral invariant of Poincar\'e-Cartan.  Given an arbitrary $2k$-dimensional phase volume
$\Omega$ (in a $2n$ dimensional space) and the Hamiltonian phase flow $\phi_t$, we have:
\[
\frac{1}{k!} \int_{\Omega} \omega^k = \frac{1}{k!} \int_{\phi_t(\Omega)} \omega^k
\]
so that the sum of the oriented $2k$-volume projections on each symplectic ``$2k$ plane'' is conserved.

\subsection{The Wirtinger-Type Integral Invariants and Volume}

Identifying $\mathbb{C}^n \simeq 
\mathbb{R}^{2n}$, the symplectic form becomes
\[
\omega = \frac{1}{2i} \sum dz_i \wedge d\overline{z}_i
\]
for any $k \in [1, N]$, consider the vectors $X_1, \ldots, X_{2k} \in \mathbb{C}^n
\simeq \mathbb{R}^{2n}$.  Wirtinger's Inequality states that the ``2k'' volume of 
the parallelpiped spanned by these $2k$ vectors is bounded by
\be
\label{sv11}
\frac{1}{k!} |\omega^k (X_1, \ldots, X_{2k})| \le \mathrm{Vol}_{2k} (X_1, \ldots, X_{2k})
\ee
We make the following two observations.  First, it is clear that $|\omega^k|$ is an integral invariant of the Hamiltonian
flow $\phi_t$, so that, given any $2k$-volume $\Omega$, we have:
\[
\frac{1}{k!} \int_{\Omega} |\omega^k| = \frac{1}{k!} \int_{\phi_t(\Omega)} |\omega^k|
\]
Moreover,
\[
\frac{1}{k!} \int_{\Omega} |\omega^k| \le \mathrm{Vol}_{2k}(\Omega)
\]
so that this integral invariant represents a minimum ($2k$) volume that the body $\Omega$ may obtain.  For the case
$k=n$, the volume of $\Omega$ is a constant which equals this invariant quantity (Liouville's Theorem).  

\subsection{Parasymplectic $2k$-Volumes}

In direct analogy with our discussion of parasymplectic surfaces, we define parasymplectic $2k$-volumes as follows:

\begin{definition} A \textbf{parasymplectic $2k$-volume}, or parametrically symplectic $2k$-volume, is one that admits
a parameterization whose paramteterization map is a symplectic one.
\end{definition}

An example of a parasympletic volume is the following.  Take any $2k$-dimensional volume $\Omega$ that is parallel to $k$ of the symplectic planes,
i.e. a region defined by:
\[
\Omega = \langle p_1, q_1, \ldots, p_k, q_k, c_{k+1}, d_{k+1}, \ldots, c_n, d_n \rangle
\]
where the variables $p_1, q_1, \ldots, p_k, q_k$ vary over some region of $\mathbb{R}^{2k}$ and $c_{k+1}, d_{k+1}, \ldots c_n, d_n$ are constants.
Now let $\phi_t$ be the Hamiltonian phase flow.  The $2k$-phase volumes $\phi_t(\Omega) \subset \mathbb{R}^{2n}$ are a one-parameter family
of parasymplectic $2k$-volumes.  

\begin{theorem} \label{svtheorem01} \textbf{[Volume Expansion of Parasymplectic $2k$-Volumes]} Let $\Omega \subset \mathbb{R}^{2k}$ be the parameterization of
a volume $\phi(\Omega) \subset \mathbb{R}^{2n}$ in a symplectic phase space whose parameterization map $\phi$ is a symplectic one.
Then 
\[
\mathrm{Vol}_{2k}(\Omega) \le \mathrm{Vol}_{2k}(\phi(\Omega))
\]
\end{theorem}

\begin{corollary} The $2k$-volume of any parasymplectic $2k$-volume is at least as large as the volume of its symplectic parameterization.
\end{corollary}
\begin{proof}
To prove our theorem, we only need prove it for a differential volume element.  The generalization follows via a simple integration argument.

Let $\phi:\mathbb{R}^{2k} \ra \mathbb{R}^{2n}$ be a symplectic parameterization of a surface. Let
\[
\left\{ \parsh{}{u_1}, \parsh{}{v_1}, \ldots, \parsh{}{u_k}, \parsh{}{v_k} \right\}
\]
be a basis of $\mathbb{R}^{2k}$ and let
\[
X_i = \phi_* \left( \parsh{}{u_i} \right) \qquad \mbox{and} \qquad Y_i = \phi_* \left( \parsh{}{v_i} \right)
\]
be the push forwards of the basis vectors in the parameterization space.  Applying Wirtinger's Inequality \eqref{sv11},
we have that
\[
\frac{1}{k!} \omega \left( \phi_* \left( \parsh{}{u_1} \right), \cdots, \phi_* \left( \parsh{}{v_k} \right) \right) \qquad 
\]
\beas
&=& \frac{1}{k!} \phi^*(\omega)\left( \parsh{}{u_1}, \ldots, \parsh{}{v_k} \right) \\
&=& \frac{1}{k!} \omega \left( \parsh{}{u_1}, \ldots, \parsh{}{v_k} \right) \\
&=& du_1 \wedge \cdots \wedge dv_k \left( \parsh{}{u_1}, \ldots, \parsh{}{v_k} \right) \\
&=& 1 \le \mathrm{Vol}_{2k} (X_1, \ldots Y_k)
\eeas
Hence any $2k$-dimensional volume measure must be nondecreasing under such a map.
\end{proof}

\subsection{The Volume Expansion Factor}

In this section we provide a practical approach to determining the volume and the integral invariants
of $2k$ subvolumes.  We will consider the volume $\Omega \subset \mathbb{R}^{2k}$ to be the 
parameterization volume of a $2k$-volume in the symplectic space $\mathbb{R}^{2n}$, with parameterization map 
\[
\phi: (\Omega \subset \mathbb{R}^{2k}) \ra (\phi(\Omega) \subset \mathbb{R}^{2n})
\]
The Jacobian of the parameterization is the $2n \times 2k$ matrix given by
\[
L = d\phi
\]
We will be interested in computing the total volume of $\Omega$, the sum of its oriented symplectic projections (i.e. 
integral invariant of Poincar\'e-Cartan), and its minimum obtainable volume (i.e. the integral invariant of the Wirtinger type).

In terms of the parameterization coordinates $\langle u_1, v_1, \ldots, u_k, v_k \rangle \in \mathbb{R}^{2k}$, the following metric is induced
on the surface:
\[
\mathbf{g} = L^T \cdot L
\]
We thus recognize the determinant of the metric $g = \det \mathbf{g}$ as the Gram determinant of the Jacobian matrix $L$:
\[
g = \mathfrak{G}(L) = \det(L^T \cdot L)
\]
so that 
\[
\mathrm{Vol}_{2k}(\phi(\Omega)) = \int_{\Omega} \sqrt{|g|} d\Omega
\]
In practical terms, the Gramian of the Jacobian can be identified with the volume expansion factor:
\[
\nu_{2k}(d\Omega; \phi) = \frac{\mathrm{Vol}_{2k}(\phi(d\Omega))}{\mathrm{Vol}_{2k}(d\Omega)} = \sqrt{\mathfrak{G}(L)}
\]
where $\nu_{2k}(d\Omega; \phi)$ is the local $2k$-volume expansion factor of $d\Omega$ under the mapping $\phi$.

\section{Local Collapse of Phase Space}

\subsection{Volume Expansion and the Local Collapse of Phase Space}

The setting for this subsection will be the evolution of a differential neighborhood surrounding a Hamiltonian trajectory
through phase space.  Consider the Hamiltonian flow:
\[ 
\phi_t: \langle p_1, q_1, \ldots, p_n, q_n \rangle \ra 
 \langle P_1, Q_1, \ldots, P_n, Q_n \rangle 
\]
Now consider a differential $2n$-``cube'' $\Omega$ situated at the initial point
$x = \langle p_1, q_1, \ldots p_n, q_n \rangle$, whose faces are parallel with the symplectic planes.
Let $\Upsilon \subset \Omega$ be a $2k$-dimensional subset that is parallel with $k$ of the symplectic planes,
and let $\Upsilon' \subset \Omega$ be a $2n-2k$ dimensional subset that is parallel with the remaining $n-k$ symplectic planes,
such that $\Omega$ is a direct sum:
\[
\Omega = \Upsilon \oplus \Upsilon'
\]
and, therefore
\[
\mathrm{Vol}_{2n}(\Omega) = \mathrm{Vol}_{2k}(\Upsilon) \cdot \mathrm{Vol}_{2n-2k}(\Upsilon')
\]
The Hamiltonian flow now takes $\phi_t: x \ra X$ along with its differential neighborhood.  We define
\beas
\overline{\Omega} &=& \phi_t(\Omega) \\
\overline{\Upsilon} &=& \phi_t(\Upsilon) \\
\overline{\Upsilon'} &=& \phi_t(\Upsilon')
\eeas
We now define the angle $\beta$ via the relation:
\[
\mathrm{Vol}_{2n} (\overline{\Omega}) = \mathrm{Vol}_{2k}(\overline{\Upsilon}) \cdot \mathrm{Vol}_{2n-2k}(\overline{\Upsilon'})
\sin \beta
\]
so that $\beta$ is the angle between the subspaces $\overline{\Upsilon}$ and $\overline{\Upsilon'}$.  
By Liouville's Theorem, we have:
\[
\mathrm{Vol}_{2n}(\overline{\Omega}) = \mathrm{Vol}_{2n}(\Omega)
\]
so that:
\[
1 = \nu_{2k}(\Upsilon; \phi_t) \nu_{2n-2k}(\Upsilon'; \phi_t) \sin \beta
\]
But by Theorem \ref{svtheorem01}, we have:
\beas
\nu_{2k}(\Upsilon; \phi_t) & \ge & 1 \\
\nu_{2n-2k}(\Upsilon'; \phi_t) & \ge & 1 
\eeas
We conclude that the greater the volume expansion of these lower dimensional differential ``slices'' $\Upsilon$ and $\Upsilon'$,
the greater the inward collapse of their respective subspaces towards each other.  In chaos theory, where 
$\mathrm{Vol}_{2k}(\phi_t(\Upsilon))$ is growing at an exponential rate, we see that $\beta$ is correspondingly decaying at an
exponential rate.  Thus chaos (for Hamiltonian systems) necessarily implies the collapse of the phase space along certain directions. 

\subsection{The Symplectic Eigenskeleton}

In this section we expose a special basis associated with any linear(ized) symplectomorphism that resists collapse.  We shall refer
to the symplectomorphism as $\phi_t$, keeping the dynamical setting (i.e. $\phi_t$ is the phase flow of a Hamiltonian system) in mind.

\begin{theorem}[The Symplectic Eigenskeleton]
Consider a symplectomorphism $\phi_t:M \ra M$ that takes the initial point $x_0$ to $\phi_t(x_0) = x$.  Let $\Phi:T_{x_0}M \ra T_x M$
be the State Transition Matrix (STM) of the mapping.  Let $\Psi = \Phi^T \cdot \Phi$ and let $\{{\xi}_1, {\eta}_1, \ldots,
{\xi}_N, {\eta}_N \}$ be the orthonormal eigenbasis of $\Psi$.  Then the following are true:
\begin{enumerate}
\item There is an interdependency amongst the vectors of $\Psi$.  The eigenvectors occur in pairs, where the $\{{\eta}_i \}_{i=1}^N$
can be taken to be
\[
{\eta}_i = J \cdot {\xi}_i
\]
where the associated eigenvalue of ${\eta}_i$ is $\la_i^{-1}$ if $\la_i$ is the eigenvalue associated with ${\xi}_i$.
\item The linear transformation $T$ that takes the standard basis to the eigenbasis of $\Psi$,
\[ T: \{ \hat{p}_1, \hat{{q}}_1,
\ldots, \hat{p}_N, \hat{q}_N \} \ra \{{\xi}_1, {\eta}_1, \ldots,
{\xi}_N, {\eta}_N \},
\]
is symplectic.  Moreover, the couples $\{ {\xi}_i, {\eta}_i \}_{i=1}^N$ make
symplectic pairs.
\item The vectors $\{ \Phi \cdot {\xi}_1, \Phi \cdot {\eta}_1, \ldots, \Phi \cdot {\xi}_N, \Phi \cdot
{\eta}_N \}$ are orthogonal.  Moreover, 
\[ 
|| \Phi \cdot \xi_i || = \sqrt{|\la_i|} \qquad \mbox{and} \qquad ||\Phi \cdot \eta_i|| = \sqrt{|\la_i^{-1}|}
\]
\item If a $2k$-dimensional symplectic subvolume $\Upsilon$ is initially parallel to $k$ of the eigenskelton planes, then the linearized transformation
$\Phi$ preserves its volume, i.e $\mathrm{Vol}_{2k}(\Upsilon) = \mathrm{Vol}_{2k}(\Phi(\Upsilon))$.
\end{enumerate}
We call the symplectic eigenbasis of the matrix $\Psi = \Phi^T \cdot \Phi$ the \textbf{symplectic eigenskeleton} of the 
transformation $\phi_t$, as it is a property structure of the transformation which resists collapse over a 
discrete time $t$.
\end{theorem}

\begin{proof}
\begin{enumerate}
\item Consider the $i$-th eigenvector $\xi_i$ of $\Psi$ with eigenvalue $\la_i$:
\[
\Psi \cdot \xi_i = \la_i \xi_i
\]
Taking the transpose of this equation, right-multiplying by $J\cdot \Psi$, and then recognizing
the identity $\Psi^T \cdot J \cdot \Psi = J$, we see that
\[
\xi_i^T \cdot J = \la_i \xi_i^T \cdot J \cdot \Psi
\]
Taking the transpose once more and multiplying by $-1$ (whilst noting $\Psi^T = \Psi$ and $J^T = -J$) we have
\[
\Psi \cdot (J \cdot \xi_i) = \frac{1}{\la_i} J \cdot \xi_i
\]
Hence, the vectors $\eta_i = J \cdot \xi_i$ are also eigenvectors of $\Psi$, with eigenvalues $\la_i^{-1}$.
\item We define the matrices $\Xi$ and $\mathrm{N}$ as follows:
\[
\Xi = \left[ \begin{array}{ccc} 
\vert &  & | \\
\xi_1 & \cdots & \xi_N \\
| &  & | \end{array} \right] \mbox{ and }
\mathrm{N} = \left[ \begin{array}{ccc} 
\vert &  & | \\
\eta_1 & \cdots & \eta_N \\
| &  & | \end{array} \right] 
\]
where $N = J\cdot \Xi$ from Part 1.  The transformation matrix $T$ can be represented as:
\[
T = \left[ \begin{array}{cc} | & | \\
\Xi & N \\
| & | \end{array} \right] 
\]
We have temporarily reordered our representation of the basis, so that 
\[
J = \left[ \begin{array}{cc}
\mathbb{O}_N & -\mathbb{I}_N \\
\mathbb{I}_N & \mathbb{O}_N
\end{array} \right]
\]
Noting again that $N = J\cdot \Xi$, one easily sees:
\[
T^T \cdot J \cdot T = \left[ \begin{array}{cc}
\Xi^T \cdot N & -\Xi^T \cdot \Xi \\
N^T \cdot N & -N^T \cdot \Xi \end{array} \right]
\]
Due to the orthonormality of the eigenvectors (i.e. $\Xi^T \cdot N = N^T \cdot \Xi = \mathbb{O}_N$ and $\Xi^T \cdot \Xi = N^T \cdot N = \mathbb{I}_N$), 
this expression reduces to:
\[
T^T \cdot J \cdot T = J
\]
and hence the matrix $T$ is symplectic.
\item Renaming the eigenvectors of $\Psi$ as $\{v_i\}_{i=1}^{2N} = \{\xi_i, \eta_i \}_{i=1}^N$, we have
from the orthonormality of the eigenbasis of $\Psi$:
\[
v_i \cdot v_j = \delta_{ij}
\]
But 
\[
(\Phi \cdot v_i) \cdot (\Phi \cdot v_j) = v_i^T \cdot \Phi^T \cdot \Phi \cdot v_j = \la_j \delta_{ij}
\]
Hence the vectors $\{ \Phi \cdot v_i \}_{i=1}^{2N}$ are also orthogonal.  Moreover:
\[
|| \Phi \cdot \xi_i || = \sqrt{ \xi_i^T \cdot \Phi^T \cdot \Phi \cdot \xi } = \sqrt{\la_i}
\]
Similarly
\[ 
|| \Phi \cdot \eta_i || = \frac{1}{\sqrt{\la_i}}
\]
\item If $\Upsilon$ is a 2-dimensional area element spanned by $\xi_1$ and $\eta_1 $, then the area expansion is
\[
\frac{\mathrm{Vol}_2(\Phi(\Upsilon))}{\mathrm{Vol}_2(\Upsilon)} = (\Phi \cdot \xi_1) \cdot (\Phi \cdot \eta_1)
= \sqrt{\la_1} \frac{1}{\sqrt{\la_1}} = 1
\]
But since this is true of any area element initially parallel to one of the symplectic eigenskeleton planes, 
and the symplectic eigenskeleton $2k$ volumes are simply direct sums of these area elements,
the result follows.
\end{enumerate}
\end{proof}

\section{Applications to Control}

We view this paper as a theoretical paper which studies some of the  
fundamental constraints in the propagation of volumes rather than  
trajectories in dynamical and control systems. This idea has already  
been advocated in viability theory and in some robust control design  
papers, see e.g. Mayne \cite{mayne}. In the final subsection below we  
present some future directions for using the theory presented here in the  
context of control. The examples below, though they do not utilize the  
full breadth of the theoretical developments presented in the paper, were  
chosen to illustrate some key ideas regarding propagation of surfaces  
and uncertainties in the control theory setting.

\subsection{The Kinematic Heisenberg System}

The Heisenberg System is a classical underactuated kinematic control problem with nonholonomic constraints,
see Bloch \cite{bloch}, Brockett \cite{brockett}.  The configuration manifold for the system is given by $Q = \mathbb{R}^3$,
with coordinates $q = \langle x, y, z \rangle$.  Motion is constrained by the relation $\dot z = y \dot x - x \dot y$.
Supposing we have controls $u, \ v$ over the $x$ and $y$ velocities, the kinematic control system can be written:
\bea
\dot x &=& u \nonumber \\
\dot y &=& v \label{sv12} \\
\dot z &=& yu-xv \nonumber
\eea
Suppose the initial state of the system is given to be within the two-dimensional uncertainty distribution
$\sigma(0) = \{ \langle x, y, z \rangle \in \mathbb{R}^3 : x \in [-1, 1], y \in [-1, 1], z = 0 \}$, 
and we wish to determine an open loop control
law that maneuvers the particle
to the point $\langle 0, 0, 1 \rangle$ during the time interval $t \in [0, 1]$ in some optimal sense.  

Let
$\sigma(t)$ be the time evolution of $\sigma(0)$ to time $t$.
We can parameterize the surface $\sigma(t)$ by the initial data $(X, Y) \in [-1, 1] \times [-1, 1]$,
so that, at time $t$, the surface is given parametrically by
$\langle x(X,Y;t), y(X,Y;t), z(X,Y;t) \rangle$.
The distance from an arbitrary point on the final surface $\sigma(1)$ to the target point is $\sqrt{x(X,Y;1)^2 + y(X,Y;1)^2 +
(1-z(X,Y;1))^2}$.  The dynamics \eqref{sv12} depend upon the choice of control $\langle u(t), v(t) \rangle$. 
We thus pose the following control problem:

\textbf{Problem:} Choose $\langle u(t), v(t) \rangle$ so as to minimize:
\be
\label{sv20}
\int \!\!\! \int_{\sigma(1)} \sqrt{x^2 + y^2 + (1-z)^2} d \mathcal{A}
\ee

\begin{figure}[t]
\begin{center}
\includegraphics[width=3.25in]{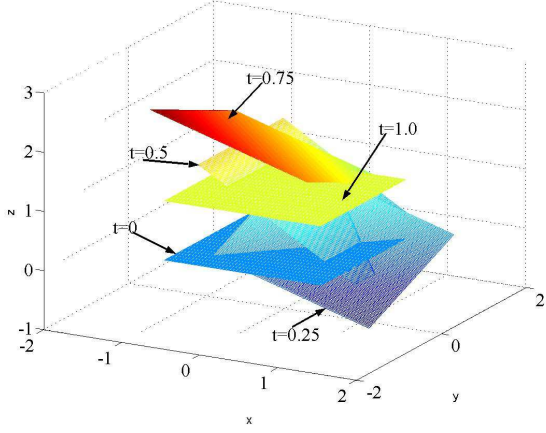}
\end{center}
\caption{Optimal Uncertainty Maneuver}
\label{svfigawesome}
\end{figure}

\textbf{Solution:}  First we need to compute the determinant of the surface metric.
The State Transition Matrix (STM) dynamics are given by:
\[ 
\dot \Phi = \left[ \begin{array}{ccc}
0 & 0 & 0 \\
0 & 0 & 0 \\
-v & u & 0 \end{array} \right] \cdot \Phi, \qquad \Phi(0) = \left[ \begin{array}{ccc}
1 & 0 & 0 \\
0 & 1 & 0 \\
0 & 0 & 1 \end{array} \right]
\]
which can immediately be integrated to yield:
\[
\Phi(t) = \left[ \begin{array}{ccc}
1 & 0 & 0 \\
0 & 1 & 0 \\
-(y(t)-y(0)) & (x(t)-x(0)) & 1 \end{array} \right]
\]
The 
metric determinant of the surface at time $t$ is given by the Gram 
determinant of the first two columns of $\Phi(t)$:
\be
\label{sv13}
g(X,Y;t) = 1 + (x(t) - X)^2 + (y(t) - Y)^2
\ee
Thus \eqref{sv20} is equivalent to:
\[
\int_{-1}^1 \int_{-1}^1 \sqrt{x^2 + y^2 + (1-z)^2} \sqrt{g(X,Y;1)} \ dX \ dY
\]
where $\langle x, y, z \rangle = \langle x(X,Y;1), y(X,Y;1), z(X,Y;1) \rangle$.
Without loss of generality, let us instead minimize:
\be
\label{sv14}
f = \int_{-1}^1 \int_{-1}^1 \left[x(1)^2 + y(1)^2 + (1-z(1))^2\right] g(1) \  dX  dY
\ee
Define now:
\[
\mu = \int_0^t u(\tau) \ d\tau,  \ \nu = \int_0^t v(\tau) \ d \tau, \ 
\alpha = \int_0^t (\nu u - \mu v) \ d \tau
\]
so that the solution to \eqref{sv12} can be expressed as:
\beas
x(t) &=& X + \mu(t) \\
y(t) &=& Y + \nu(t) \\
z(t) &=& Y \mu(t) - X \nu(t) + \alpha(t)
\eeas
This exposes the dependence of $x, \ y,$ and $z$ on the initial conditions $X$ and $Y$.  Inserting into the surface
metric \eqref{sv13}, one can explicitly integrate \eqref{sv14} to find:
\[
f = \frac 4 3 \left( 1 + \mu^2 + \nu^2 \right) \left( 4 \mu^2 + 4 \nu^2 + 3 \alpha^2 - 6 \alpha + 5 \right)
\]
This function has a global minimum at $\mu = 0, \ \nu = 0, \ \alpha = 1$.  Any control law that satisfies:
\[
\int_0^1 \!\! u(t) dt = 0,  \ \int_0^1 \!\! v(t) dt = 0,  \ \int_0^1 (\nu(t) u(t) - \mu(t) v(t)) dt = 1
\]
will leave the final uncertainty distribution as close to the target point as possible, in the above sense.
Notice that a physical interpretation of the vector quantity $\langle \mu(t), \nu(t), \alpha(t) \rangle$
is that it is the position vector of the point on the surface that was initially at $\langle 0, 0, 0 \rangle$.  
Thus, any control law that leaves the centroid of the surface at the target point will automatically minimize
\eqref{sv20}. One such trajectory, given in Bloch \cite{bloch} using $\dot y(0) = 0$, is:
\beas
\mu(t) &=& \frac{1}{\sqrt{2\pi}} \sin(2\pi t) \\
\nu(t) &=& \frac{1}{\sqrt{2\pi}} (1 - \cos(2\pi t)) \\
\alpha(t) &=& t(1 - \sin(2\pi t))
\eeas
The uncertainty surface at various time snapshots for the control law $u(t) = \dot \mu(t), \ v(t) = \dot \nu(t)$
is given in Fig. \ref{svfigawesome}.

\subsection{The Falling, Rolling Disc}

\begin{figure}[t]
\begin{center}
\psfrag{ed}{$\e_d$}
\psfrag{dphi}{$\dot{\boldsymbol{\phi}}$}
\psfrag{dtheta}{$\dot{\boldsymbol{\theta}}$}
\psfrag{dpsi}{$\dot{\boldsymbol{\psi}}$}
\psfrag{etheta}{$\e_\theta$}
\psfrag{epsi}{$\e_\psi$}
\psfrag{P}{$P$}
\psfrag{C}{$C$}
\psfrag{r}{$r$}
\psfrag{x}{$x$} \psfrag{y}{$y$}
\psfrag{z}{$z$} \psfrag{phi}{$\phi$}
\psfrag{psi}{$\psi$}
\psfrag{theta}{$\theta$}
\includegraphics[width=2.7in]{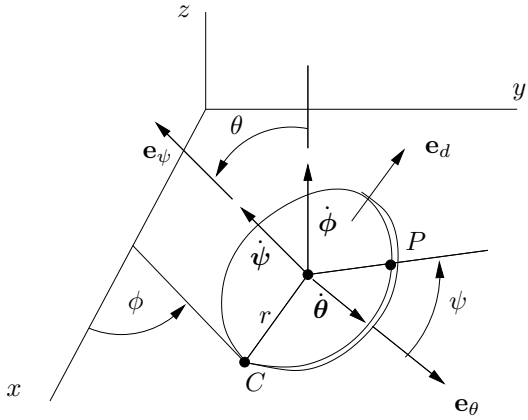}
\end{center}
\caption{Euler Angles of the Falling Rolling Disc}
\label{bhoptcon/figa06}
\end{figure}

Consider the falling rolling disc of radius $r=1$, Fig. \ref{bhoptcon/figa06}, whose configuration
is described by the contact point $(x,y)$ and the Classical Euler angles $(\phi, \theta, \psi)$.
Suppose we have direct control over the body-axis angular velocities $u = \dot \phi \sin \theta,
\ v = \dot \theta, \ w = \dot \phi \cos \theta + \dot \psi$, and suppose the system is subject to nonholonomic
constraints $\dot x + \dot \psi \cos \phi = 0$ and $\dot y +  \dot \psi \sin \phi = 0$.
The dynamics is given by the system:
\beas
\dot x &=& u \cot \theta \cos \phi  - w \cos \phi \\
\dot y &=& u \cot \theta \sin \phi - w \sin \phi \\
\dot \phi &=& u \csc \theta \\
\dot \theta &=& v \\
\dot \psi &=& -u \cot \theta  + w
\eeas
Using the notation $\dot q = f(q,u)$, the coefficient matrix in the STM dynamics equation is:
\[
\parsh{f}{q} = \left[ \begin{array}{ccccc}
0 & 0 & -(u \cot \theta - w)\sin \phi & u \csc^2 \theta \cos \phi & 0 \\
0 & 0 & (u \cot \theta - w) \cos \phi & u \csc^2 \theta \sin \phi & 0 \\
0 & 0 & 0 & - u\cot \theta \csc \theta & 0 \\
0 & 0 & 0 & 0 & 0 \\
0 & 0 & 0 & -u \csc^2 \theta  & 0
\end{array} \right]
\]
Since the STM $\Phi(t)$ is initially the identity, we find:
\[
\Phi(t) = \left[ \begin{array}{ccccc}
1 & 0 & A & C & 0 \\
0 & 1 & B & D & 0 \\
0 & 0 & 1 & E & 0 \\
0 & 0 & 0 & 1 & 0 \\
0 & 0 & 0 & F & 1
\end{array} \right]
\]
where
\beas
A &=& \int_0^t -(u \cot \theta - w) \sin \phi \ dt \\
B &=& \int_0^t (u \cot \theta - w) \cos \phi \ dt \\
C &=& \int_0^t \left( -(u \cot \theta - w) \sin \phi E + u \csc^2 \theta \cos \phi \right) \ dt \\
D &=& \int_0^t \left( (u \cot \theta - w) \cos \phi E + u \csc^2 \theta \sin \phi \right) \ dt \\
E &=& \int_0^t - u \cot \theta \csc \theta \ dt \\
F &=& \int - \csc^2 \theta \ dt
\eeas
Supposing there is initial uncertainty in the $\theta$ and $\phi$ components, it is the third and fourth columns of the STM
that will be crucial in determining the uncertainty evolution.  Suppose further our desire is that the projection of the 
final uncertainty onto the $x,y$ planes has zero area; i.e. at worst there is a one-dimensional uncertainty in the contact
point position.  Then we wish at time $t=1$, that $A(1)D(1) - B(1)C(1) = 0$.  To achieve this, one may use any control
law with that satisfies the relation 
\[
u \cot \theta - w = 0
\]
Such a control law will leave $A(t) \equiv 0, \ B(t) \equiv 0$, so that the uncertainty projection onto the $x-y$ plane
has zero area for all time.

\subsection{Future Directions in Control}

As the same symplectic constraints apply to the evolution of an optimal control system's states
and co-states, these results also have an implication for the stability and robustness
of an optimal feedback control law.  This aspect of the study can be reduced to two fundamental approaches,
the implication of initial value distributions on the subsequent evolution of a trajectory in the neighborhood 
of the true optimal trajectory, and how the symplectic invariants manifest themselves in the solutions of two-point
boundary value problems.

First, how do uncertainties in the initial state or in the initial application of the control map
to the target conditions?  As is well known, by definition an explicit optimal feedback control
law is asymptotically stable when restricted to the state variables.  However, as the necessary conditions
from which the feedback control law can arise form a Hamiltonian system, this implies that the co-states
are unstable and should diverge.  This becomes an issue if the state is not perfectly determined or
if the control function is not exact but only lies in a neighborhood of the true optimal control,
and should lead to instabilities arising in the state variables of the system.  These relationships
can be studied using integral invariants and symplectic capacities to determine the robustness of the specific
optimal control laws by studying how the phase volume surrounding them maps under the necessary conditions.
Of special interest will be the identification of the maximum and minimum uncertainty growth directions.
%squeezing directions from
%the symplectic width.

Second, given an optimal control feedback law (i.e., given the solution to the Hamilton-Jacobi-Bellman equation),
how do simultaneous uncertainties in both the initial state and target state affect the distribution of the adjoints,
and what structure may lie within these distributions that arise from the Hamiltonian formulation of the necessary
conditions?  Applying the Hamilton Principle Function approach, which provides an explicit solution to the 
two-point boundary value problem and which is directly, analytically related to the optimal control, we implicitly
define an initial set of optimal controls that will lead to a proscribed final region in the neighborhood of the nominal
target state.  This defines for us 
an open set of controls, within which lie optimal trajectories
that all achieve the final state to within some desired, 
and proscribed, accuracy.  Depending on the size and 
the distribution
of the initial uncertainties and the tolerable final uncertainties, we can identify a symplectic width which 
should provide explicit ranges in the set of initial controls that will lead to a guaranteed, optimal
arrival in the vicinity of the final state.  Such a development can provide additional insight into the 
robustness of optimal controls and how gracefully they will degrade when we allow for finite miss distances
for the target state.

There is also a clear identification between distributions in phase space and probabilistic interpretations
of the state of a system.  Thus, our research also has a direct bearing on predicted uncertainties
in a dynamical system after being mapped in time, and will define for us an absolute minimum region within
which the uncertainty of the system can be isolated.

\section{Conclusion}

We showed how the expansion of subvolumes in the local neighborhood of a nominal trajectory
leads to the local collapse of the supporting phase space.  Moreover, we produced a preferred
basis, the symplectic eigenskeleton, which resists collapse and returns uncertainty distributions
that are initially parallel to the basis to their minimal uncertainty state at a fixed final time.

%%%%%%%%%%%%%%%%%%%%%%%%%%%%%%%%%%%%%%%%%%%%%%%%%%%%%%%%%%%%%%%%%%%%%%%%%%%%%%%%
\section{ACKNOWLEDGMENTS}

We would like to thank Professor Mario Bonk, who introduced us to Wirtinger's Inequality
and who provided the initial proof of the Area Expansion of Parasymplectic Surfaces Theorem.

This research has been supported by National Science Foundation grants 
CMS-0408542 and DMS-604307.

%%%%%%%%%%%%%%%%%%%%%%%%%%%%%%%%%%%%%%%%%%%%%%%%%%%%%%%%%%%%%%%%%%%%%%%%%%%%%%%%

\end{document}